\documentclass{svproc}
\usepackage{graphicx}
\usepackage{latexsym}
 \usepackage{tabularx,ragged2e,booktabs}%,caption}
\usepackage{color}
\usepackage{float}
\usepackage{amsfonts,bbm,amsmath,amssymb}%amsthm
\usepackage{algpseudocode}
\usepackage[utf8]{inputenc}
\usepackage[T1]{fontenc}
\usepackage{enumitem}
\newenvironment{Table}
   {\par\bigskip\noindent\minipage{\columnwidth}\centering}
   {\endminipage\par\bigskip}

\begin{document}
\mainmatter
  \title{An Overview of a Break Assignment Problem Considering Area Coverage}
\author{Marin Lujak\inst{1} \and \'Alvaro Garcia S\'anchez\inst{2} \and Miguel Ortega Mier\inst{2}
\and Holger Billhardt\inst{3}}
\authorrunning{Marin Lujak et al.}
\tocauthor{Marin Lujak, \'Alvaro Garcia S\'anchez, Miguel Ortega Mier, and Holger Billhardt}
\institute{IMT Lille Douai, Univ. Lille, Unit\'e de Recherche Informatique Automatique, F-59000 Lille,
France,  \email{marin.lujak@imt-lille-douai.fr},\\
 \and
 Department of Organization Engineering, Business Administration and Statistics, Universidad
 Politécnica de Madrid (UPM), 28006 Madrid, Spain,\\
 \email{alvaro.garcia@upm.es}, \email{miguel.ortega.mier@upm.es}\\%%,
 \and
 University Rey Juan Carlos, calle Tulip\'an, s/n, M\'ostoles, Madrid, Spain,\\
 \email{holger.billhardt@urjc.es}}
\maketitle
  %\date{\today}  %\today is replaced with the current date
\begin{abstract}
Prolonged focused work periods decrease efficiency with related decline of attention and
performance. Therefore, emergency fleet break scheduling    should consider both area coverage by
idle vehicles (related to the fleet’s target arrival time to incidents) and vehicle crews' service
requirements for breaks to avoid fatigue. In this paper, we propose a break assignment problem
considering area coverage (BAPCAC) addressing this issue. Moreover, we formulate a mathematical
model for the BAPCAC problem.  Based on available historical spatio-temporal incident data and
service requirements, the BAPCAC model can be used not only to dimension the size of the fleet at
the tactical level, but also to decide upon the strategies related with break scheduling. Moreover, the model could be used to compute (suboptimal) locations for idle vehicles in each time period and arrange vehicles' crews' work breaks considering given break and coverage constraints.
 \end{abstract}
 \begingroup
\section{Introduction}\label{Introduction}
Usually, in emergency fleets (such as police and ambulance fleets), the relation between the dynamic spatio-temporal area coverage requirements and the number of breaks and their duration is not taken into account for break assignment and the breaks are regulated by (static) labor rules. Thus,  it may happen that the fleet is undersized and recurrently falls behind its performance expectations.

Even though  a vehicle crew having a break can be assigned to assist an incident, its arrival time is delayed due to the preparation time necessary to stop the break  and start the mission. This additional delay (often in terms of minutes) can mean the difference between life and death both in the case of urgent out-of-hospital patients and rapidly evolving disasters.

To guarantee efficient and timely emergency  assistance, a minimum number of stand-by vehicles with in-vehicle emergency crews must always be available  to cover each area of interest, i.e., to assist a probable incident within a predefined maximum delay in the arrival.
This can be done by dynamically transferring stand-by vehicles to  (time- and area-dependent) locations that maximize the coverage of a region of interest in each time period. These locations can be a set of strategically positioned depots, parking lots, terminals, garages, and similar.

In the emergency fleet context with temporarily and geographically varying workload, if the effective number and duration of previous work breaks  has not been taken into account in the  assignment of vehicles to incidents, it may also happen that a vehicle crew has worked during extended periods of time without any break or with  interrupted breaks that were insufficiently long for effective crew's rest.
%
%Hence, break scheduling performed  before a workshift,   discriminates crews as it does not take into account the dynamics of crews' workload as the day unfolds and related legal requirements nor it has any notion of equity among the fleet's crews nor satisfaction with work break dynamics.
Such a  disparity between a too high workload and insufficient breaks leads to increased fatigue that can have detrimental effects on emergency assistance. %The reason is that prolonged focused work periods decrease efficiency as the brain uses up oxygen and glucose.
The consequence is the feeling of being drained with related deterioration of attention and performance that may be hazardous for emergency assistance.
Hence, strategic fleet dimensioning performed before the fleet's deployment should consider expected  workload dynamics to guarantee frequent and sufficiently long time for rest for the vehicles' crews.
%Moreover, the break assignment of emergency fleet's crews should be dynamically updated through the deployment considering workload dynamics as the day unrolls and should guarantee frequent and sufficiently long time for rest. This holds especially  in 24-hour emergency shifts where eating and sleeping patterns are altered.

In this paper, we formulate this problem  under the name of the  break assignment problem  considering area coverage (BAPCAC). Based on the historical  data, the BAPCAC problem consists in determining the vehicles with in-vehicle crews necessary for optimal fleet's deployment in each time period  of a  workshift, and their breaks' type, number, and durations  based on the  legal constraints, and  the forecasted workload dynamics in terms of the area coverage  by the fleet. By the area coverage, we consider that incidents appearing in the region of interest may be potentially assisted  within a predefined  target delay in the arrival. Moreover, the BAPCAC model relocates vehicles at each time period to the locations that best cover the forecasted incident demand.
The objective  of the  BAPCAC model is an efficient dimensioning and an effective choice of break coordination strategy in emergency vehicle fleets based on historical data while increasing both the quality of service of emergency fleets as well as  the well-being of the vehicles' crew members while reducing their absenteeism due to fatigue. %Moreover,  the BAPCAC model  considers preemptive breaks, i.e.,  vehicles at break cover neighboring area considering the preparation delay time.

The paper is organized as follows. In Section \ref{RelatedWork}, we describe the state-of-the-art practice in break scheduling and introduce   the break assignment problem considering area coverage (BAPCAC) and its mathematical model. We conclude the paper with a discussion of the model in Section \ref{Discussion}.

\section{BAPCAC: Preliminaries and the model} \label{RelatedWork}
Breaks in emergency fleets are usually considered preemptive and  in the case of an  incident, the closest available vehicle having the required equipment is called to assist the incident independently of the past crew's break and workload dynamics and if it is momentarily in a break or not. This approach may induce a delay in arrival of up to several minutes, which in critical cases can be hazardous.

A review of the literature on personnel scheduling problems  can be found in \cite{van2013personnel}.
Beer et al. \cite{beer2010ai} address a complex real-world break-scheduling problem for supervisory personnel  and present a scheduling system that can help professional planners create high-quality shift plans. 
Similarly,  Di Gaspero et al.  \cite{di2013automated}  also consider the problem of scheduling breaks that fulfill different constraints about their placement and lengths.  

Mesquita et al. \cite{Mesquita2011} present  different models and algorithms for the integrated vehicle and crew scheduling problem, where crews can be assigned to different vehicles. In their algorithms, they consider several complicating constraints corresponding to workload regulations for crews.
Defraeye and Van Nieuwenhuyse \cite{defraeye2016staffing}  provide a state-of-the-art overview of research in the period 1991-2013 on personnel staffing and scheduling approaches for systems with a stochastic demand with a time-varying rate.
Rekik et al. in \cite{rekik2010implicit} consider a break scheduling problem that includes different forms of flexibility in terms of shift starting times, break lengths and break placement. 
In the  BAPCAC problem, we consider the $F|V|W$ break model. Here,  $F$ stands for the fractionable breaks with an upper and lower bound on the number of breaks that can be assigned to each vehicle crew; $V$ stands for variable  break lengths with an aggregated break length given, which  may be split into sub-breaks for each vehicle crew. $W$  stands for the workstretch duration that defines a lower and upper bound on the number of consecutive periods  of uninterrupted work before and after each sub-break can start, see, e.g., \cite{Kiermaier2014}. We also assume that the plan of \textit{vehicle crew shifts} that make a \textit{work-day shift} is given a priori.

We assume that incidents appear, w.l.o.g., in a 2D square environment $Env=[0, l]^2 \subset \mathbf{R} ^2$ of side length $l >0$. Environment $Env$ is divided in mutually exclusive and nonempty rectangular unit areas $e \subset Env$ of equal size such that $e_i \cap e_j= \emptyset$ for $i\neq j$ and $\bigcup_{e \in Env} e = Env$. Vehicle agents are positioned and may move only  within the subset $\mathcal{J}  \subseteq Env$ of unit areas.

Let $\mathcal{T}$ be a set of consecutive time periods $\{1, \ldots, |\mathcal{T}|\}$, representing a vehicle crew shift, where  $|\mathcal{T}|$ is the cardinality of $\mathcal{T}$. For each time period  $\tau \in \mathcal{T}$, $D_{e\tau}$  represents  incident density or demand at area $e \in Env$. 
$\mathcal{I} \subseteq Env$ is a subset of unit areas $i \in \mathcal{I}$  with positive incident density.  
Then, each vehicle $a \in A$ can have one of  mutually exclusive vehicle states: $\{idle, on-break, occupied\}$: \emph{idle} -  a vehicle is waiting for a new incident assistance;
\emph{on-break} - currently on a break. 
An idle vehicle is considered to be covering unit area $i \subset \mathcal{I}$ if it is located in any unit area from which it can arrive to unit area $i$  within a predefined and given target arrival time $\Delta t$, while for the vehicles on-break, the arrival time is increased by an additional preparation delay.

Let $B=\{b_s, b_l\}$ be a set of  short and long breaks respectively.
The minimum number of short and long breaks is related to the duration of the shift, labor rules, and legal requirements.  We track accumulated workload of each vehicle crew since the last break  at each time period $\tau \in \mathcal{T}$  to comply with a given maximum allowed work time  $MAX_b^w$ before the start of a break $b \in B$. 
Each break has a minimal duration $\Delta_b^{MIN}$ and a maximal duration $\Delta_b^{MAX}$, that specify the required amount of break time based on legal requirements.

Since the problem of area coverage is interconnected with the break assignment problem in the emergency fleet context, we  approach it by a unique mathematical programming model in which the decisions on the number of vehicles necessary for the coverage, their locations, identities, and the break assignments are decided all at once for the whole workshift.

\textit{Sets and indices}

\begin{Table}
  \begin{tabular}{ll}
$\mathcal{T}$ & time horizon; a set of time periods in a work shift; $\tau \in \mathcal{T}$\\
$A$& set of agents $a \in A$ representing $|A|$ capacitated vehicles \\
	&with vehicle crews  \\
$B$ & set of break types for example,  $\{b_{s}, b_l\}$, where\\
	& $b_{s}$ stands for short breaks, whereas $b_l$ stands for long breaks\\
%$Env=[0,l]^2$ & 2D square environment divided into unit areas \\
$\mathcal{I}$ & set of unit areas $i \in \mathcal{I}$ with positive incident density\\
$\mathcal{J}$ & set of unit areas $j \in \mathcal{J}$, apt for vehicle locations \\
  \end{tabular}
    %  \caption{Sets in the problem definition}  \label{Table_parameters}
  \end{Table}
\textit{Parameters}

\begin{Table}
  \begin{tabular}{ll}
   $D_{i\tau}$ & estimated demand and, thus, required minimum proportion of vehicles \\
				& for unit area $i \in \mathcal{I}$ at time $\tau\in \mathcal{T}$\\
   %$\Delta t$ & predefined target arrival time \\ %This is not used in the formulation
	 $\Delta^{MIN}_b$& minimal duration of break  $b \in B$\\%, where  $\Delta_B \leq U_b-L_b$\\
	 $\Delta^{MAX}_b$& maximal duration of break $b \in B$\\
  	 $IM_{jj'}$& valued 1 if an agent can move from $j\in \mathcal{J}$ to $j'\in \mathcal{J}$  based on \\
				& the maximum transport time between 2 consecutive  periods; \\
& 0 otherwise\\
   $N_{ij}$ &  1 if an idle vehicle in $j \in \mathcal{J}$ can assist an incident at $i \in \mathcal{I}$ \\
				& within a given target arrival time, 0 otherwise\\
   $\bar{N}_{ij}$ & 1 if a vehicle on a break located at $j \in \mathcal{J}$ can assist an incident \\
   & at $i \in \mathcal{I}$ within a predefined target arrival time, 0 otherwise\\% $\Delta t$\\
   $MAX^w_b$ &  maximum allowed work time before assigning break $b \in B$ \\
		&usual this time will be 3 hours for short and 4-6 hours for long breaks\\
  \end{tabular}
  \end{Table}
\textit{Decision variables}

\begin{Table}
  \begin{tabular}{ll}
 $x_{a\tau j}$ & valued 1 if agent a at time $\tau \in \mathcal{T}$ is located \\
			& at unit area $j \in \mathcal{J}$; 0 otherwise\\
 $y_{a\tau}$ & binary break assignment variable valued $1$ if agent $a \in A$ \\
			& is assigned a break at time $\tau \in \mathcal{T}$; 0 otherwise\\
 $z_{a\tau i}$ & real nonnegative coverage assignment variable representing \\
 			& the part of density at incident unit area $i \in \mathcal{I}$ assigned to\\
			& \textit{idle}  agent $a \in A$ at time $\tau \in \mathcal{T}$, where  $0\leq z_{a\tau i} \leq 1$\\
 $\bar{z}_{a \tau i}$& real nonnegative coverage assignment variable representing \\
 			& the part of density at incident unit area $i \in \mathcal{I}$ assigned to \\
 			& \textit{on-break} agent $a \in A$  at time $\tau \in \mathcal{T}$, where  $0\leq \bar{z}_{a\tau i} \leq 1$ \\
$\delta_{i\tau}$& real nonnegative density of area $i \in \mathcal{I}$ uncovered at time $\tau \in \mathcal{T}$\\
$\alpha_{ab\tau} $ & binary variable valued 1 if agent $a \in \mathcal{A}$ starts a break of type $b \in B$ \\
& at $\tau \in \mathcal{T}$; 0 otherwise\\
   \end{tabular}
  \end{Table}

\indent (BAPCAC):
\begin{equation}\label{ObjectiveFunction}
\min  \;w\cdot \sum_{i \in \mathcal{I}, \tau\in \mathcal{T}} \delta_{i\tau}  + (1-w)\cdot \sum_{a\in A, \tau\in \mathcal{T}} \big(1-y_{a\tau}\big),
\end{equation}

such that:
\begin{equation}\label{BAPCAC Coverage}
\sum_{a \in A}\big( z_{a\tau i} + \bar{z}_{a \tau i} \big) = D_{i\tau} - \delta_{i\tau}, \forall i \in \mathcal{I}, \tau \in \mathcal{T}
\end{equation}

\begin{equation}\label{BAPCAC_Individ_CoverageLin1}
z_{a \tau i} \leq  \sum_{j \in \mathcal{J}|N_{ij}=1}D_{i\tau}x_{a \tau j}, \forall a \in A, \tau \in \mathcal{T}, i \in \mathcal{I}
\end{equation}

\begin{equation}\label{BAPCAC_Individ_CoverageLin2}
\bar{z}_{a \tau i} \leq  \sum_{j \in \mathcal{J}|\bar{N}_{ij}=1}D_{i\tau}x_{a \tau j}, \forall a \in A, \tau \in \mathcal{T}, i \in \mathcal{I}
\end{equation}

\begin{equation}\label{BAPCAC_Individ_CoverageLin3}
z_{a \tau i} \leq  D_{i\tau}(1-y_{a,\tau}), \forall a \in A, \tau \in \mathcal{T}, i \in \mathcal{I}
\end{equation}

\begin{equation}\label{BAPCAC_Individ_CoverageLin4}
\bar{z}_{a \tau i} \leq  D_{i\tau}y_{a,\tau}, \forall a \in A, \tau \in \mathcal{T}, i \in \mathcal{I}
\end{equation}

\begin{equation}\label{BAPCAC_Coverage_1Agent}
\sum_{i \in \mathcal{I}}z_{a \tau i} \leq 1, \forall \tau \in \mathcal{T}, a \in A
\end{equation}

\begin{equation}\label{BAPCAC_Coverage_OnBreak}
\sum_{i \in \mathcal{I}}\bar{z}_{a \tau i} \leq 1, \forall \tau \in \mathcal{T}, a \in A
\end{equation}

\begin{equation}\label{BAPCAC_Position_in_Area}
\sum_{j \in \mathcal{J}} x_{a \tau j} = 1, \forall \tau \in \mathcal{T}, a \in A
\end{equation}

\begin{equation}\label{BAPCAC_Max_Move_Btw_Areas}
x_{a \tau j} +  x_{a (\tau+1) j'} \leq 1 + IM_{jj'} , \forall a \in A, \tau =\{1, |\mathcal{T}|-1\}, j,j' \in \mathcal{J}
\end{equation}

\begin{equation}\label{CLSCP-min-duration}
\sum_{\tau'=\tau}^{\tau+\Delta^{MIN}_b} y_{a\tau'} \geq \Delta^{MIN}_b \alpha_{ab\tau}, \,\forall a \in A,\, b \in B,\,\tau =\{1, |\mathcal{T}|-\Delta^{MIN}_b\}
\end{equation}

%Constr.12
\begin{equation}\label{CLSCP-max-duration}
y_{a\,(\tau+\Delta^{MAX}_b+1)}+ \alpha_{ab\tau}  \leq 1 , \,\forall a \in A,\, b \in B,\,\tau =\{1, |\mathcal{T}|-\Delta^{MAX}_b-1\}
\end{equation}

\begin{equation}\label{CLSCP-short-break-starts}
\sum_{\tau'=\tau}^{\tau+MAX^w_b}\alpha_{ab\tau'}= 1 , \, \forall a \in A,\, b \in B, \tau \in \{1, \mathcal{T}-MAX_b^w\}
\end{equation}

\begin{equation}\label{OverlappingBreaks}
\sum_{\tau'=\tau}^{\tau+MAX_{b_s}^w}\alpha_{a b_s \tau'} +
\sum_{\tau'=\tau}^{\tau+MAX_{b_l}^w}\alpha_{a b_s \tau'} \leq 2,  \forall a \in A,   \tau \in \{1,  \mathcal{T}-MAX_{b_l}^w\}
\end{equation}

\begin{equation}\label{CLSCP NonnegFinal2}
 z_{a\tau i},\bar{z}_{a\tau i}, \delta_{i \tau}\geq 0,   x_{a\tau j}, y_{a \tau}, \alpha_{a b \tau} \in \{0,1\}, \forall i \in \mathcal{I}, j \in \mathcal{J},  \tau \in \mathcal{T}, a \in A
\end{equation}

\section{Discussion}\label{Discussion}
In (\ref{ObjectiveFunction}), $w$  is the  weight assigned to the coverage mean  such that $0\leq w \leq 1$. Constraints (\ref{BAPCAC Coverage}) relate to the coverage of each unit area $i$ by \textit{idle} vehicles and  vehicles \textit{on-break}.
Here, minimizing the  arrival time to an incident  has  priority over the performance of a break in its totality.
If we exclude $\bar{z}_{a \tau i}$ from the above and other constraints, our model turns into a non-preemptive break model.

Constraints (\ref{BAPCAC_Individ_CoverageLin1}) and (\ref{BAPCAC_Individ_CoverageLin2}) limit the coverage of each idle agent and  an agent on-break  to at most density $D_{i\tau}$ of  unit area $i$ if it is positioned within the travel time defined by adjacency matrix $N_{ij}$ and $\bar{N}_{ij}$,  respectively.
These constraints ensure that part $z_{a \tau i}$ and $\bar{z}_{a \tau i}$   of the coverage of area $i$  by an idle agent and an agent on-break $a$ respectively,  is at most the sum of the densities of the areas $j \in \mathcal{J}$ within its reach at time $\tau \in \mathcal{T}$.
Furthermore, constraints (\ref{BAPCAC_Individ_CoverageLin3}) and (\ref{BAPCAC_Individ_CoverageLin4}) limit idle coverage and the coverage on-break only to idle  agents and the agents at break, respectively.
Moreover, constraints (\ref{BAPCAC_Coverage_1Agent}) %related to the formulation in (\ref{BAPCAC_Break_Coverage})
guarantee that the incident density covered by each idle agent  sums up to at most  1.
Similarly, constraints (\ref{BAPCAC_Coverage_OnBreak}) capacitate  the coverage  for each agent on-break.

Constraints (\ref{BAPCAC_Position_in_Area}) assign to each agent $a$ at each time period $\tau$ a unique location $j \in \mathcal{J}$, by means of $x_{a \tau j}$, while  constraints (\ref{BAPCAC_Max_Move_Btw_Areas}) relate agents' positions $x_{a \tau j}$ in two consecutive time periods $t$, $t+1$ based on the maximum distance travelled  established in $IM_{jj'}$.
Moreover, minimal and maximal break durations are imposed by constraints (\ref{CLSCP-min-duration}) and (\ref{CLSCP-max-duration}).
Constraints (\ref{CLSCP-short-break-starts}) assure that short and long breaks start at latest after a maximum uninterrupted work time  $MAX^w_{b}$ of consecutive time periods without a break.
Additionally, (\ref{OverlappingBreaks}) specifies that long and short breaks are coordinated within the maximum allowed working times. Finally, constraints (\ref{CLSCP NonnegFinal2}) are non-negativity constraints on the decision variables. 
\paragraph{Acknowledgements} This work has been partially supported  by an STSM Grant from COST Action
TD1409 ``Mathematics for Industry Network (MI-NET)''.
%\section*{References}
\bibliographystyle{splncs03}
\bibliography{lunch_breaks}
\endgroup
\end{document}